\documentclass{amsart}
\usepackage{graphicx}        
\RequirePackage{hyperref}   
\providecommand{\href}[2]{#2}
\begin{document}
%
%
%
%
%
{ 
 
\title{\bf Spatial Filtering for Reduced Order Modeling}
\author{L.~C.~Berselli \and D.~Wells \and X.~Xie \and T.~Iliescu}
\address[L.~C.~Berselli]{Dip. di Matematica, Univ. di Pisa,
  Pisa, Italy} \email {luigi.carlo.berselli@unipi.it}

\address[D.~Wells] {Dept. of Math. Sci., Rensselaer Polytechnic Inst.,
  Troy, NY, US} 
\email{wellsd2@rpi.edu} 
\address[X.~ Xie and T.~Iliescu]{ Dept. of Mathematics,
  Virginia Tech., Blacksburg, VA, US}
\email{iliescu@vt.edu}
%
\maketitle
%
%
%
%
%
\providecommand{\ff}{\vec{f}}
\providecommand{\uu}{\vec{u}}
\providecommand{\va}{\vec{a}}
\providecommand{\vv}{\vec{v}}
\providecommand{\ww}{\vec{w}}
\providecommand{\FF}{\mathcal{F}}
\providecommand{\podvec}{\vec{\varphi}}
\providecommand{\eg}{e.g., }
\providecommand{\ie}{i.e., }
\providecommand{\OO}{\mathcal{O}}
\section{Introduction}
Spatial filtering has been central in the development of large eddy
simulation reduced order models
(LES-ROMs)~\cite{wang2012proper,xie2017data,xie2017approximate} and
regularized reduced order models
(Reg-ROMs)~\cite{sabetghadam2012alpha,wells2017evolve} for efficient
and relatively accurate numerical simulation of convection-dominated
fluid flows.  In this paper, we perform a numerical investigation of
spatial filtering.  To this end, we consider one of the simplest
Reg-ROMs, the Leray ROM
(L-ROM)~\cite{sabetghadam2012alpha,wells2017evolve}, which uses ROM
spatial filtering to smooth the flow variables and decrease the amount
of energy aliased to the lower index ROM basis functions.  We also
propose a new form of ROM differential
filter~\cite{sabetghadam2012alpha,wells2017evolve} and use it as a
spatial filter for the L-ROM.  We investigate the performance of this
new form of ROM differential filter in the numerical simulation of a
flow past a circular cylinder at a Reynolds number $Re=760$.
\section{Reduced Order Modeling}
For the Navier-Stokes equations (NSE), the standard reduced order
model (ROM) is constructed as follows: (i) choose modes $\{ \podvec_1,
\ldots, \podvec_d \}$, which represent the recurrent spatial
structures of the given flow; (ii) choose the dominant modes $\{
\podvec_1, \ldots, \podvec_r \}$, $r \leq d$, as basis functions for
the ROM; (iii) use a Galerkin truncation $\uu_r = \sum_{j=1}^{r} a_j
\, \podvec_j$; (iv) replace $\uu$ with $\uu_r$ in the NSE; (iii) use a
Galerkin projection of NSE($\uu_r$) onto the ROM space $X^r :=
\mbox{span} \{ \podvec_1, \ldots, \podvec_r \}$ to obtain a
low-dimensional dynamical system, which represents the ROM:
\begin{eqnarray}
	\dot{\va}
	= A \, \va
	+ \va^{\top} \, B \, \va \, ,
	\label{eqn:g-rom}
\end{eqnarray}
where $\va$ is the vector of unknown ROM coefficients and $A, B$ are
ROM operators; (iv) in an offline stage, compute the ROM operators;
and (v) in an online stage, repeatedly use the ROM (for various
parameter settings and/or longer time intervals).
\section{ROM Differential Filter}
\label{sec:filtering:pod-differential-filter}
The ROM differential filter is based on the classic Helmholtz filter
that has been used to great success in LES for turbulent
flows~\cite{germano1986differential}. Let \(\delta\) be the radius of
the differential filter. Then, for a given velocity field \(\uu_r \in
X^r\), the filtered flow field \(\FF(\uu_r) \in X^f\), where \(X^f\)
is a yet to be specified space of filtered ROM functions, is defined
as the solution to the Helmholtz problem
\begin{equation}
  {\rm Find}\  \FF(\uu_r) \in X^f \ {\rm such\ that}\ 
  \label{eq:filtering:pod-differential-filter}
  \left((I - \delta^2 \Delta) \FF(\uu_r), \vv\right)
  = (\uu_r, \vv),
  \ {\rm for\  all}\  \vv \in X^f.
\end{equation}
We consider two different versions for the choice of the range of the
ROM differential filter \(X^f\):
\smallskip

\noindent\textbf{The FE Version.} This version corresponds to \(X^f =
X^h\), where $X^h$ is the finite element (FE) space: we seek the FE
representation of \(\FF(\uu)\) and work in the full discrete space
when calculating the filtered ROM vectors. The FE
representation of \(\FF(\uu)\) suffices in applications because we use
it to assemble the components of the ROM before time evolution: put
another way, since filtering is a linear procedure, it only has to be
done once and not in every ROM time step, \eg for FE mass
and stiffness matrices \(M\) and \(S\) we have that, modulo boundary
condition terms,
\begin{equation}
  a_j (M + \delta^2 S) \FF\left(\podvec_j\right)
  = a_j M \podvec_j
  \Rightarrow
  (M + \delta^2 S) \sum_{j = 1}^r a_j \FF(\podvec_j)
  = M \sum_{j = 1}^r a_j \podvec_j.
  \label{eq:filtering:discrete-pod-differential-filter}
\end{equation}
Hence, applying the differential filter to each proper orthogonal
decomposition (POD) basis vector \(\podvec_j\), results in
\(\FF(\podvec_j) \notin X^r\).
Due to the properties of the differential filter (see
Fig.~\ref{figure:filtering:filtered-pod-vectors}), these new ROM
functions will correspond to longer length scales and contain less
energy.  \smallskip

\noindent\textbf{The ROM Version.} Alternatively, we can pick \(X^f =
X^r\), \ie the ROM differential filter simply corresponds to an \(r
\times r\) Helmholtz problem. 
\begin{equation}
  (M_r + \delta^2 S_r) \FF\left(\va\right)
  = M_r \va \, ,
  \label{eq:rom-df}
\end{equation}
where $M_r$ and $S_r$ and the ROM mass and stiffness matrices,
respectively, and $\va$ and $\FF\left(\va\right)$ are the POD
coefficient vectors of \(\podvec_j\) and \(\FF(\podvec_j)\),
respectively. 
Here, unlike in the FE version, the range of the Helmholtz filter is
\(X^r\), so filtered solutions retain the weakly divergence free
property.  \smallskip

\noindent\textbf{Properties.}
\begin{figure}[t]
\centering
\includegraphics[width=1.0\textwidth]
{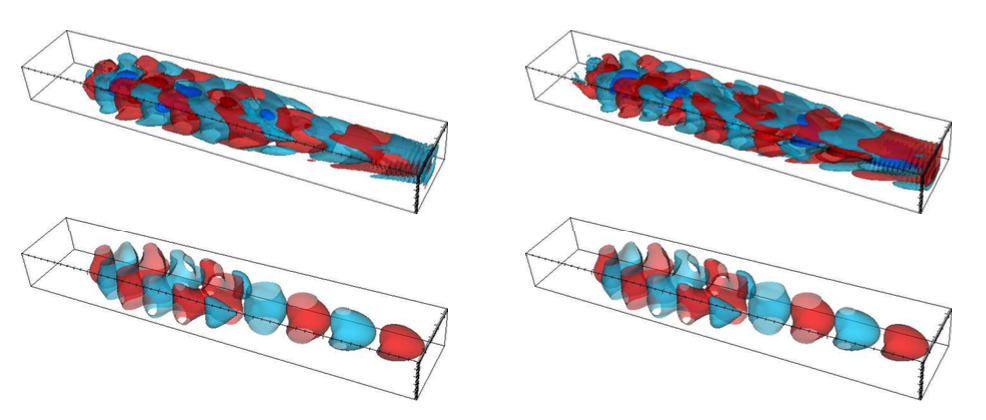}
\caption{Contour plots of \(y\) velocity of the first and fifth POD
  vectors from a 3D flow past a cylinder computation. We use the FE
  version of the ROM differential filter with $\delta = 0.5$. The
  unfiltered POD vectors are on top and the filtered are on the
  bottom. The choice $\delta = 0.5$ is too large for practical
  purposes, but demonstrates that filtering both removes kinetic
  energy (the isosurfaces are smaller) and enlarges the scales of
  motion (\eg the first POD vector goes from twelve structures to just
  nine.)}
        \label{figure:filtering:filtered-pod-vectors}
\end{figure}
Both versions of the ROM differential
filter~(\ref{eq:filtering:pod-differential-filter}) share several
appealing properties~\cite{BIL05}. They act as spatial filters, since
they eliminate the small scales (\ie high frequencies) from the input.
Indeed, the ROM differential
filter~(\ref{eq:filtering:pod-differential-filter}) uses an elliptic
operator to smooth the input variable.  They also have a low
computational overhead.  For efficiency, the algorithmic complexity of
any additional filters should be dominated by the \(\OO(r^3)\) cost in
evaluating the nonlinearity. The ROM version is equivalent to solving
an \(r \times r\) linear system; since the matrix only depends on the
POD basis, it may be factorized and repeatedly solved for a cost of
\(\OO(r^2)\), which is also dominated by the cost of the nonlinearity.
The FE version requires solving large FE linear systems, but these
linear systems are solved in the offline stage; thus, the online
computational cost of the FE version is negligible.  Finally, we
emphasize that the ROM differential filter uses an \emph{explicit}
length scale \(\delta\) to filter the ROM solution vector.  This is
contrast to other types of spatial filtering, e.g., the ROM
projection, which do not employ an explicit length scale.
\section{Leray ROM}
\label{sec:regrom:leray-model}
Jean Leray attempted to solve the NSEs in his landmark 1934
paper~\cite{leray1934sur}. He was able to prove the existence of
solutions for the modified problem
\begin{equation}
  \label{eq:regrom:leray_regularized_ns_pde}   
  \ww_t = \frac{1}{Re} \Delta \ww-\FF(\ww) \cdot \nabla \ww-\nabla p,                                                         
 \end{equation}
 where $\nabla \cdot \ww =0$, and \(\FF(\ww)\) is a convolution with a
 compact support mollifier with filter radius \(\delta\), or
\begin{equation}
  \label{eq:regrom:gaussian_convolution}
  \FF(\ww) = g_\delta \star \ww.
\end{equation}
For additional discussion on the properties of different filters
see~\cite{BIL05,layton2012approximate, Sag06}. 
We approximate the convolution with the differential filter
\begin{equation}
  \label{eq:regrom:differential_filter}
  \FF(\ww) = (\delta^2 \Delta + 1)^{-1} \ww.
\end{equation}
In turbulence modeling, Leray's model is the basis for a class of
stabilization methods called the
Leray-\(\alpha\) regularization models~\cite{layton2012approximate}. Leray's
key observation was that the nonlinear term is the most problematic as
it serves to transfer energy from resolved to unresolved
scales. 

The Leray model has been recently extended to the ROM
setting~\cite{sabetghadam2012alpha,wells2017evolve}.  The resulting
\emph{Leray-ROM (L-ROM)} can be written as
\begin{equation}
  \label{eq:regrom:filtered_centered_pde}
  (\ww_r)_t = \frac{1}{Re} \Delta \ww_r-\FF(\ww_r) \cdot \nabla \ww_r-\nabla p,                                                         
\end{equation}
which is the same as the Galerkin ROM up to the filtering of the
advective term in the nonlinearity.

\section{Numerical Results}
We consider the flow past a cylinder problem with parabolic Dirichlet
inflow conditions, no-slip boundary conditions on the walls of the
domain, and zero tangential flow at the outflow. %
We compute snapshots by running the \texttt{deal.II}
\cite{bangerth2016deal} step-35 tutorial program for \(t \in [0, 500]\). 
We use a kinematic viscosity value of $1/100$,  a circular
cylinder with diameter of $1$, and parabolic inflow boundary 
conditions with a maximum velocity of $7.6$; this results in a 
Reynolds number $Re=760$.  
We calibrate the filter radius \(\delta\) by choosing a value for
\(\delta\) that gives the L-ROM the same mean kinetic energy as the
original numerical simulation.  Calibrating the ROM to this filter
radius also improves accuracy in some structural properties: this
amount of filtering removes enough kinetic energy that the phase
portrait connecting the coefficients in the ROM on the first and
second POD basis functions are close to the values obtained by
projecting the snapshots onto the POD basis over the same time
interval.

Fig.~\ref{fig:regrom:ns-leray-df-l2-norms} displays the time evolution
of the \(L^2\) norm of the solutions of the L-ROM and DNS for \(r =
6\) and \(r = 20\).  Fig.~\ref{fig:regrom:ns-leray-df-l2-norms} shows
that, for the optimal \(\delta\) value, the L-ROM-DF accurately
reproduces the average, but not the amplitude of the time evolution of
the \(L^2\) norm of the DNS results for both \(r = 6\) and \(r = 20\).
\begin{figure}[!htb]
        \label{fig:regrom:ns-leray-df-l2-norms}
        \centering
\includegraphics[width=1.0\textwidth]
{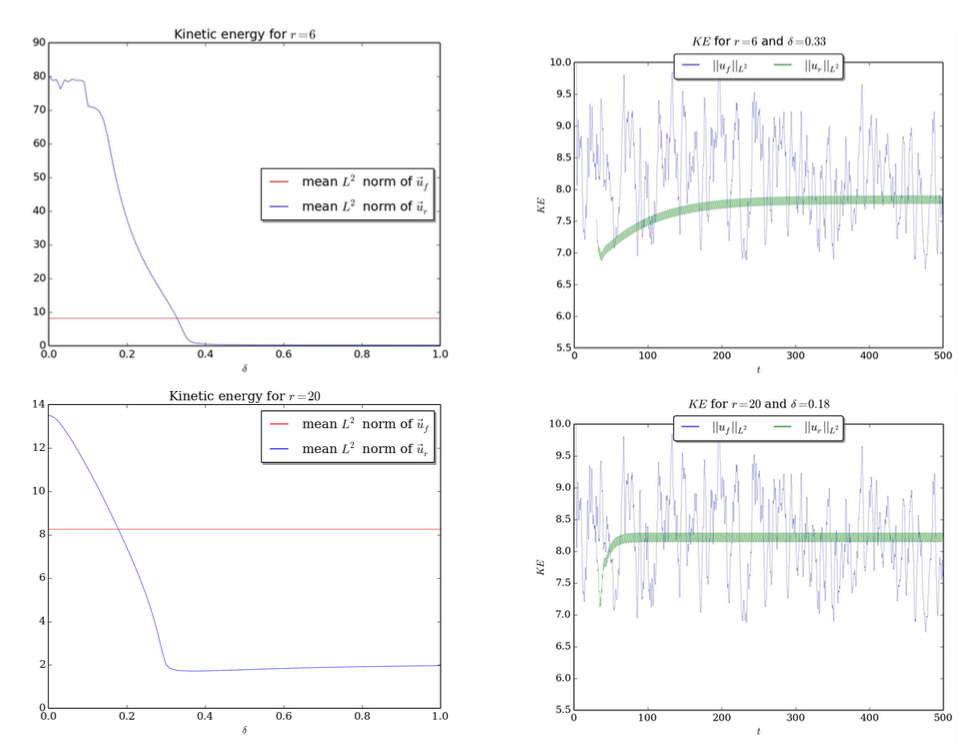}
        \caption{3D flow past a cylinder, L-ROM (green) and DNS
          (blue). Mean (left column) and time evolution (right column)
          of the $L^2$ norm of the solution; $r = 6$ (top row) and $r
          = 20$ (bottom row). The time evolution of the $L^2$ norm of
          the solution (right column) is plotted for the optimal mean
          $L^2$ norm of the solution (left column): $\delta = 0.33$
          for $r = 6$ (top row) and $\delta = 0.18$ for $r = 20$
          (bottom row).  }
    \end{figure}
    Fig.~\ref{fig:regrom:ns-leray-df-phase-0-1} displays the phase
    portraits for the first and second POD coefficients of the
    L-ROM-DF and POD projection of DNS data for \(r = 6\) and \(r =
    20\). Fig.~\ref{fig:regrom:ns-leray-df-phase-0-1} shows that, for
    the optimal \(\delta\) value, the L-ROM-DF yields moderately
    accurate results for \(r = 6\) and accurate results for \(r =20\).
    \begin{figure}[!htb]
        \label{fig:regrom:ns-leray-df-phase-0-1}
        \centering
\includegraphics[width=1.0\textwidth]
{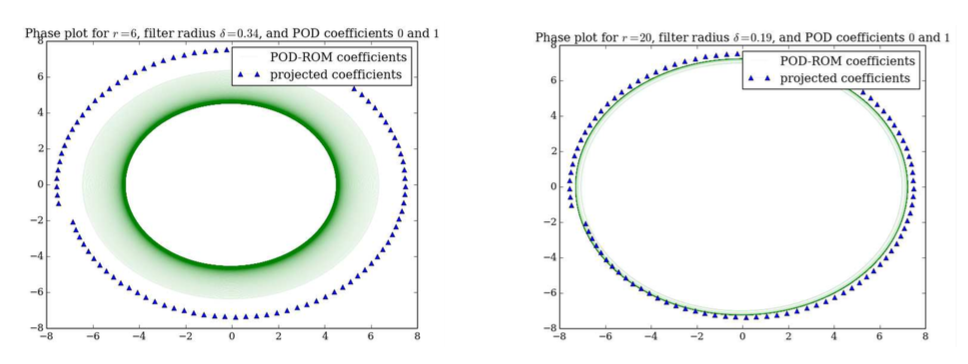}
        \caption{3D flow past a cylinder, L-ROM-DF with optimal
          $\delta$ value (green) and POD projection of DNS data
          (blue). Phase portraits for $a_1$ and $a_2$; $r = 6$ (left)
          and $r = 20$ (right).}
    \end{figure}
\section{Conclusions}
In this paper, we proposed a new type of ROM differential filter.  We
used this new filter with the L-ROM, which is one of the simplest
Reg-ROMs.  We tested this filter/ROM combination in the numerical
simulation of a flow past a circular cylinder at Reynolds number
$Re=760$ for $r=6$ and $r=20$.  The new type of ROM differential
filter yielded encouraging numerical results, which were comparable to
those for the standard type of ROM differential filter and better than
those for the ROM projection~\cite{wells2017evolve}.  We emphasize
that a major advantage of the new type of ROM differential filter over
the standard ROM differential filter is its low computational
overhead.  Indeed, since the filtering operation in the new type of
ROM differential filter is performed at a FE level (as opposed to the
ROM level, as it is generally done), the new filter is applied to each
ROM basis function in the offline stage.  In the online stage, the
computational overhead of the new type of ROM differential filter is
practically zero, since it simply amounts to using the filtered ROM
basis functions computed and stored in the offline stage.

The first results for the new type of ROM differential filter are
encouraging.  We plan to perform a thorough investigation of the new
filter, including a comparison with the standard form of the ROM
differential filter and the ROM projection, in the numerical
simulation of realistic
flows~\cite{wells2017evolve,xie2017approximate}.

}
\end{document}